\documentclass[a4paper]{article}

\usepackage{amssymb}
\usepackage{amsmath}
\usepackage{latexsym}
\usepackage[all]{xy}

\title{Generating Specials: The Zorro Algorithm}
\author{Jacob Thamsborg}
\date{\today}

\newtheorem{defn}{Definition}
\newtheorem{lem}[defn]{Lemma}
\newtheorem{prop}[defn]{Proposition}
\newtheorem{thm}[defn]{Theorem}
\newtheorem{cor}[defn]{Corollary}

\begin{document}
\maketitle

\begin{abstract}
The concept of a configuration graph associated to a primitive,
aperiodic substitution is introduced in \cite{eilers} as a
convenient graphical representation of the infinite indeterminism of
the shift space of the substitution. The main result of
\cite{eilers} is an algorithm to calculate this graph from the
substitution, in this paper we turn the tables and produce
substitutions from graphs. We do this using the Zorro algorithm, an
entirely constructive and easily applicable algorithm. In the
process we show that any configuration graph can be obtained.

The first section contains standard definitions and the definition
of configuration graphs. The second and third sections develop
theory used in the proof of the algorithm as stated in section four.
The algorithm is easily applied without knowledge of the underlying
theory. Note that section three is nothing but a copy of results
from \cite{eilers} slightly modified to suit the present needs.
\end{abstract}

\section{Preliminaries}
\subsection{Meeting Notational Needs}
Let $\mathcal{A}$ be any nonempty finite set of symbols, we call
$\mathcal{A}$ our \emph{alphabet} and its members \emph{letters}. By
$\mathcal{A}^*$ we understand the set of finite words constructed
from the letters of $\mathcal{A}$ including the empty word
$\epsilon$. Equipped with the associative composition of
concatenation, $\mathcal{A}^*$ is the free monoid over
$\mathcal{A}$. We furthermore let $\mathcal{A}^+ = \mathcal{A}^*
\backslash \{\epsilon\}$ denote the set of nonempty words, and for
any $u \in \mathcal{A}^*$ we let $|u|$ be the length of $u$, i.e.,
the number of letters of $u$. Given two words $u$ and $v$ of
$\mathcal{A}^*$ we say that $u$ is a \emph{factor} of $v$ denoted $u
\dashv v$ if there exists $w_1, w_2 \in \mathcal{A}^*$ with $w_1uw_2
= v$.

We call members of $\mathcal{A}^\mathbb{Z}$ \emph{(two sided)
sequences} over the alphabet $\mathcal{A}$. Let $x$ be some sequence
and let $i \in \mathbb{Z}$, we denote the letter at index $i$ with
$x_{[i]}$, given an additional $j \in \mathbb{Z}$ with $i \leq j$ we
let $x_{[i,j]}$ denote the word consisting of the letters from index
$i$ to index $j$, both included. We define the \emph{language} of
some sequence $x$ to be the set $\mathcal{L}(x) = \{\epsilon \} \cup
\left\{u \in \mathcal{A}^* \mid \exists i,j \in \mathbb{Z}, i \leq
j: u = x_{[i,j]}\right\}$ and call its members \emph{factors} of
$x$. We define the \emph{shift} $\sigma:\mathcal{A}^\mathbb{Z} \to
\mathcal{A}^\mathbb{Z}$ by $\left(\sigma(x)\right)_{[i]}=x_{[i+1]}$
for $x \in \mathcal{A}^\mathbb{Z}$ a sequence and $i$ ranging over $
\mathbb{Z}$. Elements of $\mathcal{A}^\mathbb{N}$ are called
\emph{one sided sequences} over $\mathcal{A}$; subscript notation
and definition of language, factors and shift apply to these as
well, only the indices range over $\mathbb{N}$ and not $\mathbb{Z}$.
Note, however, that while the shift is bijective on
$\mathcal{A}^\mathbb{Z}$ it is only surjective on
$\mathcal{A}^\mathbb{N}$.

Let $u$ be any word of $\mathcal{A}^*$ and $x$ a one sided sequence,
the concatenation $ux$ is defined the obvious way. Given a two sided
sequence $x$ and $i \in \mathbb{Z}$ we let $x_{]\infty,i]}$ and
$x_{[i,\infty[}$ denote obvious one sided sequences. Given, on the
other hand, any two single sided sequences $x$ and $y$, we define
the two sided sequence $x.y$ by letting $x.y_{[i]} = x{[-i]}$ for $i
< 0$ and $x.y_{[i]} = y_{[i+1]}$ for $i \geq 0$, i.e., by reversing
$x$ and concatenating it with $y$, letting the first letter of $y$
have index $0$. We shall extend this notation in the obvious way to
allow for finite words between the dot and the one sided sequences.
For the sake of an example, let $x$ and $y$ be one sided sequences
and let $a$ be some letter, we then have that $\sigma(x.ay) = xa.y$.

By a \emph{substitution} $\tau$ we understand a map
$\tau:\mathcal{A} \to \mathcal{A}^+$, it can be extended in the
obvious way to a map respecting concatenation $\tau:\mathcal{A}^*
\to \mathcal{A}^*$, furthermore to map single sided sequences to
single sided sequences and by specifying $\tau(x.y) =
\tau(x).\tau(y)$ for any $x,y \in \mathcal{A}^\mathbb{N}$ to map
sequences to sequences; we shall not distinguish between a
substitution and its extension. Note that for any $u \in
\mathcal{A}^*$ we have $|\tau(u)| \geq |u|$ and that for any two
substitutions $\tau_1$ and $\tau_2$ the composition $\tau_1 \tau_2$
defines a substitution as well.

\subsection{Primitivity and Aperiodicity: Pretty Interesting Substitutions}
In this subsection we introduce the concept of primitivity, the
language associated with a substitution, the shift space associated
with a substitution and finally the concept of aperiodicity. The
different properties are easily verified if one proceeds in the
order they are listed here.

\begin{defn}
A substitution $\tau$ is said to be primitive if it holds that
\[\exists n \in \mathbb{N} \forall a, b \in \mathcal{A}: b \dashv
\tau^n(a) \] and that
\[\exists a \in \mathcal{A} \forall N \in \mathbb{N} \exists n \in
\mathbb{N}: |\tau^n(a)| > N. \]
\end{defn}
Notice that the first of these properties implies the second if we
have $|\mathcal{A}| > 1$, indeed the second property does nothing
but exclude the substitution $a \mapsto a$ in a theoretically
convenient way.
\begin{prop}
Let $\tau$ be any primitive substitution. We have the following
properties:
\begin{itemize}
\item[(i)] $\exists n \in \mathbb{N} \forall a, b \in \mathcal{A} \forall i \in \mathbb{N}_0: b \dashv
\tau^{n+i}(a)$
\item[(ii)] $\forall a \in \mathcal{A} \forall N \in \mathbb{N} \exists n \in
\mathbb{N}: |\tau^n(a)| > N$
\item[(iii)] $\exists x \in \mathcal{A}^\mathbb{Z}\exists n \in \mathbb{N}:
\tau^n(x) = x$
\end{itemize}
\end{prop}

Now let $\tau$ be some substitution, we define the \emph{langauge}
of $\tau$ by
\[\mathcal{L}(\tau) = \left\{u \in \mathcal{A}^* \mid \exists a \in
\mathcal{A}\exists n \in \mathbb{N} : u \dashv \tau^n(a) \right\}.
\]

\begin{prop}
Let $\tau$ be any substitution. We have the following properties:
\begin{itemize}
\item[(i)] $\tau(\mathcal{L}(\tau)) \subseteq
\mathcal{L}(\tau)$
\item[(ii)] $\forall u, v \in \mathcal{A}^*: u \dashv v, v \in
\mathcal{L}(\tau) \Rightarrow u \in \mathcal{L}(\tau)$
\end{itemize}
If furthermore $\tau$ is primitive we get that:
\begin{itemize}
\item[(iii)] $\mathcal{A} \subseteq \mathcal{L}(\tau)$
\item[(iv)] $\forall n \in \mathbb{N}: \mathcal{L}(\tau) = \mathcal{L}(\tau^n)$
\end{itemize}
\end{prop}
Consider now the non primitive substitution:
\[ \tau_{d}: 1 \mapsto 2, \ 2 \mapsto 3, \ 3 \mapsto 3. \]
We obviously have $\mathcal{L}(\tau_d) = \{\epsilon, 2, 3\}$ and
$\mathcal{L}(\tau_d^2) = \{\epsilon, 3\}$ which demonstrates that
primitivity is a necessary condition for the two lower properties.

We furthermore define the \emph{shift space} associated with $\tau$
by
\[ X_\tau = \left\{x \in \mathcal{A}^\mathbb{Z} \mid \mathcal{L}(x)
\subseteq \mathcal{L}(\tau) \right\}. \]

\begin{prop}
Let $\tau$ be any substitution. We have the following properties:
\begin{itemize}
\item[(i)] $\sigma(X_\tau) = X_\tau$
\item[(ii)] $\tau(X_\tau) \subseteq X_\tau$
\end{itemize}
If furthermore $\tau$ is primitive we get that:
\begin{itemize}
\item[(iii)] $\forall x \in X_\tau \forall u \in \mathcal{L}(x) \exists
n \in \mathbb{N}_0 \forall i \in \mathbb{Z}: u \dashv x_{[i,i+n]}$
\item[(iv)] $\forall x \in X_\tau: \mathcal{L}(x) = \mathcal{L}(\tau)$
\item[(v)] $\forall n \in \mathbb{N}: X_\tau = X_{\tau^n}$
\item[(vi)] $X_\tau \neq \emptyset$
\end{itemize}
\end{prop}

A sequence $x \in \mathcal{A}^\mathbb{Z}$ is said to be
\emph{periodic} if there exists an $n \in \mathbb{N}$ such that for
all $i \in \mathbb{Z}$ we have $x_{[i]} = _{[i+n]}$, $n$ is called
the length of the period. Finally let $\tau$ be a primitive
substitution. We say that $\tau$ is \emph{periodic} if $X_\tau$ is
finite. This is equivalent to $\tau$ having a periodic member of
$X_\tau$ which is again equivalent to having all members of $X_\tau$
periodic. \emph{Aperiodicity} is obviously defined as the lack of
periodicity for sequences as well as substitutions.

We end this somewhat tedious subsection with a small but handy
lemma:
\begin{lem} \label{perlem}
Let $\tau$ be a primitive, aperiodic substitution and let $x \in
X_\tau$. We have that
\[\forall i,j \in \mathbb{Z}:x_{[i,\infty[} = x_{[j,\infty[}
\Leftrightarrow i = j \]
\end{lem}
The proof is an easy application of the definitions above, a
symmetrical version of the lemma also holds.

%last by the following substitution:
%\[ \tau_{d'}: 1 \mapsto 11, \ 2 \mapsto 22, \]
%which obviously has $\cdots 111.111 \cdots \in X_{\tau_{d'}}$ but
%which also has $2 \in \mathcal{L}(\tau_{d'})$ and $2 \notin
%\mathcal{L}(\cdots 111.111 \cdots)$. Again, primitivity will help us
%out here as we shall see in the next section.

\subsection{Orbit classes, specials and configuration graphs}
\label{specials}
Let $\tau$ be a primitive, aperiodic substitution.
By definition this implies that $X_\tau$ is infinite. In this
subsection we shall consider the structure of $X_\tau$, in
particular we shall present the concept of a configuration graph
associated to $\tau$ which is a convenient graphical representation
of the infinite indeterminism of $X_\tau$.

\begin{defn}
Let $\tau$ be any substitution. Let $x, y \in X_\tau$. We define the
following relations:
\begin{itemize}
\item[(i)] $x \sim_o y \Leftrightarrow \exists m \in \mathbb{Z}
\forall i \in \mathbb{Z}: x_{[i]} = y_{[i+m]}$
\item[(ii)] $x \sim_r y \Leftrightarrow \exists m \in \mathbb{Z}
\exists M \in \mathbb{Z} \forall i \geq M: x_{[i]} = y_{[i+m]}$
\item[(iii)] $x \sim_l y \Leftrightarrow \exists m \in \mathbb{Z}
\exists M \in \mathbb{Z} \forall i \leq M: x_{[i]} = y_{[i+m]}$
\end{itemize}
\end{defn}
We name these relations \emph{orbit equivalence}, \emph{right tail
equivalence} respectively \emph{left tail equivalence} and
immediately verify that they are indeed equivalence relations. The
equivalence classes under orbit equivalence are called \emph{orbit
classes} and since both right and left tail equivalence respect
orbit equivalence they define equivalence relations on the orbit
classes as well.

\begin{defn}
Let $\tau$ be any substitution. A sequence $x \in X_\tau$ is called
\emph{left special} if there exists $y \in X_\tau$ with
\[ x_{[-1]} \neq y_{[-1]} \quad x_{[0,\infty[} = y_{[0,\infty[}. \]
An orbit class $C \in X_\tau/\sim_o$ is called \emph{left special}
if there exists an orbit class $D \in X_\tau/\sim_o$ with $C \neq D$
and $C \sim_r D$.
\end{defn}
And yes, an easy application of lemma \ref{perlem} shows that if
$\tau$ is primitive and aperiodic then an orbit class is left
special if and only if it contains a left special sequence.
\emph{Right special} sequences and orbit classes are defined
symmetrically.

As mentioned in theorem 1.5 of \cite{eilers} the number of left as
well as right special orbit classes is finite but nonzero if $\tau$
is primitive and aperiodic. This makes the following definition
meaningful:
\begin{defn}
Let $\tau$ be a primitive, aperiodic substitution. The
\emph{configuration graph} is a bipartite graph defined as follows:
The set of left vertices are the equivalence classes of orbit
classes under left tail equivalence that contain a special orbit
class. The set of right vertices are defined symmetrically and each
special orbit class gives rise to an edge connecting the left and
right equivalence classes that contain it.
\end{defn}

As an example, the primitive, aperiodic substitution $1 \mapsto
121,\ 2 \mapsto 2112$ has the following configuration graph:
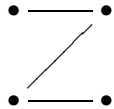
\[\xymatrix{\bullet\ar@{-}[r]&\bullet\\\bullet\ar@{-}[ur]\ar@{-}[r]&\bullet}\]
The calculation of configuration graphs is by no means a trivial
exercise, indeed an algorithm doing this is the main result of
\cite{eilers}. This algorithm is most conveniently implemented
online, see \cite{applet} for details.

\section{Generators}
\begin{defn}
Let $\tau$ be any substitution. Let $(v,u,w) \in \mathcal{A}^+
\times \mathcal{A}^+ \times \mathcal{A}^+$. We say that $(v,u,w)$ is
a generator for $\tau$ if $u \in \mathcal{L}(\tau)$ and furthermore
$\tau(u) = vuw$. We denote by $G_{\tau}$ the set of all generators
for $\tau$.
\end{defn}
Given a generator $(v,u,w)$ we shall refer to $v$, $u$ and $w$ as
the \emph{left wing}, the \emph{center} respectively the \emph{right
wing} to facilitate the language. Furthermore we shall refer to the
length of the center as the length of the generator.
\begin{defn}
Let $\tau$ be any substitution and let $(v,u,w) \in G_{\tau}$. We
define the completion of $(v,u,w)$ by $$ (v,u,w)^{*} = \cdots
\tau^2(v)\tau(v)vu.w\tau(w)\tau^2(w)\cdots$$ and note that this is a
member of $X_{\tau}$.
\end{defn}
This definition is our main justification for working with
generators: they provide a means of creating members $X_\tau$ and
they do so in a nice way as we shall see below. But before we start
completing let us first impose some structure on the set of
generators.

\begin{defn}
Let $\tau$ be any substitution and let $(v,u,aw)$ be a generator
with $v, u \in \mathcal{A}^+$, $a \in \mathcal{A}$ and $w \in
\mathcal{A}^*$. Then obviously $(v,ua,w\tau(a))$ is a generator as
well and we say it it constructed from the original by \emph{right
extension}; \emph{left extension} is defined similarly. We say that
two generators $g_1$ and $g_2$ for $\tau$ are G related (denoted by
$g_1 \sim_{G} g_2$) if there exists a generator $g_3$ such that
$g_3$ can be constructed from $g_1$ by a series of (possibly zero)
right and left extensions and $g_3$ can be constructed similarly
from $g_2$.
\end{defn}
One quickly realizes that left as well as right extensions are
deterministic, i.e., any generator can be left or right extended in
exactly one way. Furthermore, right and left extensions are
independent since they take place on different sides of the center,
so to speak, and this implies that their order can be exchanged in a
a series of mixed extensions. Summing up, the relation defined above
is transitive as well as obviously reflexive and symmetric, i.e., it
is an equivalence relation.

\begin{defn}
Let $\tau$ be any substitution. We define the basic generators to be
all generators that are not G related to any shorter generator.
\end{defn}
We shall see shortly that there is exactly one basic generator in
each equivalence class. But let us pause to consider how we would
calculate the basic generators of a substitution, this turns out to
be very easy in the case of primitive substitutions:
\begin{lem}
Let $\tau$ be any substitution and let $g = (v,aub,w)$ be any
generator of two or more letters with $v, w \in \mathcal{A}^+$, $u
\in \mathcal{A}^*$ and $a,b \in \mathcal{A}$. It is basic if and
only if $|\tau(a)| > |v|$ and $|\tau(b)| > |w|$.
\end{lem}
\emph{Proof:} Suppose one of the length inequalities fail, say,
$|\tau(a)| \leq |v|$. Then we can write $v = \tau(a)v'$ for some $v'
\in \mathcal{A}^*$ and $(v'a,ub,w)$ is a generator shorter than $g$
and obviously G related to $g$.

Now suppose both length inequalities hold. Let $n \in \mathbb{N}_0$.
We shall show by complete induction on $n$ that if $g'$ and $g''$
are two more generators and $g$ can be extended to $g''$ in a series
of $n$ extensions and $g'$ can be extended to $g''$ in another
series of extensions, then $g'$ is longer than or has the same
length as $g$. Let $m$ be the number of left extensions of the $n$
steps and let $m'$ be the number of left extensions in the steps
extending $g'$ to $g''$. If both are nonzero we can remove one left
extension from both series and still end up with a common result,
since left and right extensions commute, and afterwards apply the
inductive hypothesis. We cannot have $m = 0$ and $m' > 0$ since the
first would let the left length inequality hold for $g''$ and the
second would contradict this. This leaves us with $m \geq m'$ and
since the same arguments applies to right extensions we have
finished our inductive argument and the proof.
\hspace{\stretch{1}}$\Box$
\begin{cor}
Let $\tau$ be a primitive substitution. The following holds:
\begin{itemize}
\item[(i)] All one letter generators are basic.
\item[(ii)] Let $(v,ab,w)$ be any two letter generator with
$v, w \in \mathcal{A}^+$ and $a,b \in \mathcal{A}$. It is basic
if and only if $\tau(a) = va$ and $\tau(b) = bw$.
\item[(iii)] No generators of three or more letters are basic.
\end{itemize}
\end{cor}
Notice that the primitivity condition is necessary for part (iii)
since a non primitive substitution may have basic generators of any
length. Consider for instance the following non primitive
substitution:
\[ 0 \mapsto 01230, \ 1 \mapsto 1, \ 2 \mapsto 2, \ 3 \mapsto 30123. \]
This has the generator $(0123,0123,0123)$ which is basic by the
lemma thus contradicting the corollary.

The set of basic generators of a primitive substitution is very
easily calculated using the corollary: The one letter generators can
be read off the definition of the substitution directly; the two
letter generators in question are all those that can be constructed
from mating a one letter empty right wing "generator" with a one
letter empty left wing "generator", bearing in mind that the center
must always be in $\mathcal{L}(\tau)$. As an example consider the
following primitive substitution:
\[ 0 \mapsto 042, \ 1 \mapsto 142, \ 2 \mapsto 042, \ 3 \mapsto 043, \ 4 \mapsto 01432. \]
This has the four basic generators $(01,4,32)$, $(04,20,42)$,
$(04,21,42)$ and $(04,30,42)$ and no more, in particular
$(04,31,42)$ is not even a generator.

The following proposition justifies the basic generators as being,
in essence, all generators:
\begin{prop}
Let $\tau$ be any substitution. We then have:
\begin{itemize}
\item[$(i)$] No two different basic generators are G related.
\item[$(ii)$] Any generator is G related to a unique basic
generator.
\end{itemize}
\end{prop}
\emph{Proof:} The proof of (i) proceeds similarly to the proof of
the second part of the lemma, i.e., complete induction on the number
of steps required to extend $g$ to some generator that another basic
generator can be extended to as well. Common left extensions are
handled by the inductive hypothesis and left extensions in only one
of the extension series are contradicted by the lemma. The proof of
(ii) is immediate by induction on the length of the generator by the
definition of basic generators; the uniqueness is a spinoff from
part (i). \hspace{\stretch{1}}$\Box$

It is now time to consider how these structures on $G_\tau$ interact
with the completion of members of $G_\tau$. The following result is
a pretty one:
\begin{prop} \label{simsim}
Let $\tau$ be a primitive, aperiodic substitution and let $g_1, g_2
\in G_{\tau}$. We have that
\[ g_1^{*} \sim_o g_2^{*} \iff g_1 \sim_G g_2. \]
\end{prop}
\emph{Proof:} The arrow leading left is immediate since left and
right extension preserve completion up to orbit equivalence.

Assume now that $g_1^{*} \sim_o g_2^{*}$. Assume initially that
$g_1^{*} = g_2^{*}$. Let $n_1, n_2 \in \mathbb{N}$ be the length of
the right wing of $g_1$ respectively $g_2$. Since
\[
\sigma^{-n_1}(\tau(g_1^{*})) = g_1^{*} = g_2^{*} =
\sigma^{-n_2}(\tau(g_2^{*})) = \sigma^{-n_2}(\tau(g_1^{*}))
\]
aperiodicity ensures that $n_1 = n_2$. This immediately implies that
if $g_1$ and $g_2$ are of equal length then they are equal, and if
they are not, then the shorter can be left extended to obtain
longer. If $g_1^{*} \neq g_2^{*}$ then there must exist a $p \in
\mathbb{Z}, p \neq 0$ such that $\sigma^p(g_1^{*}) = g_2^{*}$. In
case $p > 0$ then by performing $p$ right extensions of $g_1$ we are
in the situation above. The case $p<0$ is handled by right extending
$g_2$.\hspace{\stretch{1}}$\Box$

With the construction of specials in mind, the following result is
promising:
\begin{prop} \label{simr}
Let $\tau$ be a primitive, aperiodic substitution and let $g_1, g_2
\in G_{\tau}$. We have that $g_1^{*} \sim_r g_2^{*}$ holds if and
only if there exist two generators $g'_1 \sim_G g_1$ and $g'_2
\sim_G g_2$ with identical right wings.
\end{prop}
\emph{Proof:} Assume that $g_1^{*} \sim_r g_2^{*}$ holds. If we have
the luck that $g_1 \sim_G g_2$ then by definition there exists
 a $g'$ with $g_1 \sim_G g'$ and $g_2 \sim g'$ and letting $g'_1 =
 g'$ and $g'_2 = g'$ concludes the case. If, on the other hand, $g_1 \nsim_G
 g_2$ holds then we have the existence of $p, j \in \mathbb{Z}$ such that
\[\forall i \geq j: \sigma^p(g_1^{*})_{[i]} = {g_2^{*}}_{[i]} \]
and
\[\sigma^p(g_1^{*})_{[j-1]} \neq {g_2^{*}}_{[j-1]}.\]
Assume initially that $p = 0$. If we further assume assume that $j
\leq 0$, then we can halfway duplicate the calculations from the
proof of proposition \ref{simsim}: Let $n_1, n_2 \in \mathbb{N}$ be
the length of the right wing of $g_1$ respectively $g_2$. We now
get:
\begin{eqnarray*}
\sigma^{-n_1}(\tau(g_1^{*}))_{[n_2, \infty[} & = & {g_1^{*}}_{[n_2, \infty[} \\
                             & = & {g_2^{*}}_{[n_2, \infty[} \\
                             & = & \sigma^{-n_2}(\tau(g_2^{*}))_{[n_2, \infty[} \\
                             & = & \sigma^{-n_2}(\tau(g_1^{*}))_{[n_2, \infty[}
\end{eqnarray*}
This by lemma \ref{perlem} is enough to ensure that $n_1=n_2$ which
proves that the two generators have identical right wings. Now if
$j>0$ then we perform $j$ right extensions on both generators and
proceed as above, this concludes the case $p=0$. And as above, if $p
> 0$ then we do $p$ right extensions of $g_1$, if $p<0$ then we do
$p$ right extensions of $g_2$ and in both cases proceed as in the
case $p=0$. The reverse is immediate. \hspace{\stretch{1}}$\Box$

Given two basic generators $g_1$ and $g_2$ with $g_1 \nsim_G g_2$
and suppose we'd like to know whether $g_1^{*} \sim_r g_2^{*}$. The
proposition above tells us to look for G related generators with
identical right wings, but this is not an algorithmically very
pleasant task. But the proof above shows that $g'_1$ and $g'_2$ --
if they exist at all -- can be constructed by doing nothing but
right extensions of $g_1$ respectively $g_2$. After possibly undoing
some pairwise identical right extensions we can furthermore obtain
generators with identical right wings that disagree on either their
rightmost letter of the center or the letter just before that. If
now additionally $\tau$ is regular, then this puts a maximum limit
to the length of the desired common right wing, thereby making the
test for $g_1^{*} \sim_r g_2^{*}$ a finite story. Let us list an
even simpler and most useful case:
\begin{cor} \label{simrid}
Let $\tau$ be a primitive, aperiodic, postfix free substitution and
let $g_1, g_2 \in G_{\tau}$ with $g_1 \nsim_G g_2$. We have that
$g_1^{*} \sim_r g_2^{*}$ holds if and only if the right wings of
$g_1$ and $g_2$ are identical.
\end{cor}

A final note to conclude this section: The definition of the
completion of a generator is not entirely symmetrical with respect
to the left and right wings of the generator. The given definition
has the pleasant property that right extending the generator shifts
the completion one step; we rely heavily on this in the proofs
above. On the other hand, one might fear that this would introduce
some asymmetry to completions. This, however, is not the case as
long as we stick to orbit classes. Indeed, the symmetrical versions
of both proposition \ref{simr} and corollary \ref{simrid} above
hold, this is most easily checked by shifting to opposite
substitutions.

\section{Generating specials}
\begin{defn}
Let $\tau$ be any substitution. The leftmost letter graph (the ll
graph) is defined to be the graph with the letters of $\mathcal{A}$
as vertices and with one directed edge leaving each vertex $a \in
\mathcal{A}$ arriving at the leftmost letter of $\tau(a)$. The
rightmost letter graph (the rl graph) is defined similarly.
\end{defn}

\begin{defn}
Let $\tau$ be any substitution and let $n \in \mathbb{N}$. We say
that $n$ is a left segregating number if for any two words $u, v \in
\mathcal{L}_n(\tau)$ with differing leftmost letter we have that the
length of the common prefix of $\tau(u)$ and $\tau(v)$ is less than
or equal to $min\left\{|\tau(u)|, |\tau(v)|\right\} - n$. Right
segregating numbers are defined similarly.
\end{defn}

Note that not all substitutions have a segregating numbers. Consider
for instance the following primitive, aperiodic substitution:
\[ \tau_e: a \mapsto c, \ b \mapsto c, \ c \mapsto db, \ d \mapsto ca.
\] Squaring this we get a substitution with the two generators
$(d,bca,cdb)$ and $(ca,cdb,cdb)$. This implies that for any $n \in
\mathbb{N}$ there exists $u \in A^*$ with $|u| = n-1$ and $au, bu
\in \mathcal{L}_n(\tau_e)$ which shows that $n$ cannot be a left
segregating number since we have that $\tau_e(au) = \tau_e(bu)$. On
the other hand, note that for any prefix free substitution $1$ will
do as left segregating number, similarly any postfix free
substitution has $1$ as right segregating number. We say that a
substitution is \emph{segregating} if it has both a left and a right
segregating number. As is often the case, regular substitutions
behave nicely:

\begin{prop} \label{regseg}
Let $\tau$ be any primitive, regular substitution. Then $\tau$ is
segregating.
\end{prop}

\emph{Proof:} We prove only the existence of the left segregating
number, the right case is symmetrical. Since $\tau$ is primitive
there must exist an $a \in \mathcal{A}$ with $\tau(a) > 1$. By
minimality there exists an $s \in \mathbb{N}$ such that any $u \in
\mathcal{L}_s(\tau)$ contains $a$. Now let
\[ P = \sum_{a \in \mathcal{A}} |\tau(a)|,\  Q = \max_{a \in \mathcal{A}} |\tau(a)|.\]
It now follows from theorem 1.6 in \cite{salomaa} that $s(P -
\left|\mathcal{A} \right| + Q -1)$ is a left segregating number.
\hspace{\stretch{1}}$\Box$

\begin{defn}
Let $\tau$ be any substitution with a left segregating number. Let
$n \in \mathbb{N}$ be the least such. We define the left segregating
graph (the ls graph) as follows: The vertices are all pairs of words
from $\mathcal{L}_n(\tau)$ which differ at their leftmost letter.
One directed edge leaves each vertex, if the vertex is $(u,v)$ then
the destination is obtained by removing the common prefix from
$\tau(u)$ and $\tau(v)$ and reading the leftmost $n$ letters from
each remaining word. The right segregating graph (the rs graph) is
defined similarly for a substitution with a right segregating
number.
\end{defn}

It is time for an example, consider the following primitive,
aperiodic, regular substitution:
\[ \tau_4: 0 \mapsto 10, \ 1 \mapsto 0. \]
The ll and rl graphs are as follows:
\[\xymatrix{\mathrm{ll:} & 0
\ar@/^/[r] & 1 \ar@/^/[l] & \mathrm{rl:} & 0 \ar@(ul,dl)[] & 1.
\ar@/^/[l] } \] As left segregating number 1 will do, and clearly it
is the least such. On the other hand, 2 is the least right
segregating number. Since $\mathcal{L}_1(\tau) = \{0, 1\}$ and
$\mathcal{L}_2(\tau) = \{00, 01, 10\}$ we get the following ls and
rs graphs:
\[\xymatrix{\mathrm{ls:} & (0,1) \ar@/^/[r] & (1,0) \ar@/^/[l]}
\qquad \xymatrix{\mathrm{rs:} & (00,01) \ar[dr] & (10,01) \ar@/^/[d]
\\ & (01,00) \ar[ur] & (01,10). \ar@/^/[u]} \]

We say that any of the graphs defined above are \emph{subfixed} if
for each vertex $v$, $v$ either loops to itself (i.e., the edge
leaving $v$ goes back to $v$) or the edge leaving $v$ goes to some
other vertex that loops to itself. Of the graphs in the example
above only the rl graph is subfixed. It is, however, the case that
for any segregating substitution $\tau$ there exists an $n \in
\mathbb{N}$ such that all the graphs ll, rl, ls and rs for $\tau^n$
are subfixed. To realize this, notice first that if $\tau$ is
segregating then so is any nonzero power of $\tau$. Then note that
raising the power of $\tau$ by one corresponds to extending each
edge by its immediate successor in any of the graphs above. Finally
let $m$ be the least common multiple of the length of all cycles in
all the graphs (each must have at least one cycle if $\tau$ is
primitive and aperiodic). Then raising $\tau$ to the power of any
positive multiple of $m$ ensures that all vertices that are in
cycles the original graph now loop to themselves and by choosing a
sufficiently high multiple we can make all other vertices connect to
one of these vertices. In the simple example above choosing $n=2$
will work, i.e., for $\tau_4^2$ all the graphs ll, rl, ls and rs are
subfixed. The following theorem is our main justification for this
as well as the preceding section:

\begin{thm} \label{exists}
Let $\tau$ be any primitive, aperiodic, segregating substitution
with all the graphs ll, rl, ls and rs subfixed. Then for any left or
right special sequence $u \in X_\tau$ there exists a generator $g
\in G_\tau$ such that $g^* \sim_o u$.
\end{thm}

To prove this, consider first the following lemma:
\begin{lem} \label{adopt}
Let $\tau$ be any primitive, aperiodic substitution with a right
segregating number and with the rs graph subfixed. Suppose we have
$u, v \in X_\tau$ with $u_{[0,\infty[} = v_{[0,\infty[}$ and
$u_{[-1]} \neq v_{[-1]}$. Then there exist $u', v' \in X_\tau$ with
$u'_{[0,\infty[} = v'_{[0,\infty[}$ and $u'_{[-1]} \neq v'_{[-1]}$
and
\[ u[-n,-1] = u'[-n,-1], \ v[-n,-1] = v'[-n,-1] \]
and
\[ u = \sigma^{-p}\left(\tau(u')\right), \ v = \sigma^{-p}\left(\tau(v')\right),  \]
where $n \in \mathbb{N}$ is the least right segregating number and
$p \in \mathbb{N}_0$ is the length of the common postfix of
$\tau(u'[-n,-1])$ and $\tau(v'[-n,-1])$.
\end{lem}

\emph{Proof of lemma:} By corollary 12 of \cite{durand} there exists
$x, y \in X_\tau$ with $u \sim_o \tau(x)$ and $v \sim_o \tau(y)$. By
lemma 3.1 of \cite{eilers} we get that $x \sim_r y$. But since $u
\nsim_o v$ we also have $x \nsim_o y$ and we may choose $u' \sim_o
x$ and $v' \sim_o y$ with $u'_{[0,\infty[} = v'_{[0,\infty[}$ and
$u'_{[-1]} \neq v'_{[-1]}$. Now there exists $p,q \in \mathbb{Z}$
such that $u = \sigma^{-p}(\tau(u'))$ and $v =
\sigma^{-q}(\tau(v'))$ but it follows from lemma \ref{perlem} that
$p=q$ and we can furthermore deduce that these must equal the length
of the common postfix of $\tau(u'[-n,-1])$ and $\tau(v'[-n,-1])$.
Now repeat this exercise to produce $u''$ and $v''$ with
$u''_{[0,\infty[} = v''_{[0,\infty[}$ and $u''_{[-1]} \neq
v''_{[-1]}$ and with $u'=\sigma^{-r}(\tau(u''))$
$v'=\sigma^{-r}(\tau(v''))$ where $r$ is the length of the common
postfix of $\tau(u''[-n,-1])$ and $\tau(v''[-n,-1])$. Now going from
$(u'',v'')$ to $(u',v')$ and on to $(u,v)$ makes the pair of words
at index $[-n,-1]$ change according to the rs graph and since this
is subfixed we have that $u[-n,-1] = u'[-n,-1]$ and $v[-n,-1] =
v'[-n,-1]$ as desired. \hspace{\stretch{1}}$\Box$

\emph{Proof of theorem:} We assume that $u$ is left special, the
right case is, as is often the case, symmetrical. By definition
there must exist $v \in X_\tau$ with $u_{[0,\infty[} =
v_{[0,\infty[}$ and $u_{[-1]} \neq v_{[-1]}$. Now let $n \in
\mathbb{N}$ be the least right segregating number, let $p \in
\mathbb{N}_0$ be the length of the common postfix of
$\tau(u[-n,-1])$ and $\tau(v[-n,-1])$ and let $r \in \mathbb{N}_0$
be $|\tau(u[-n,-1])| - p - n$. Now suppose both $p$ and $r$ are
nonzero. Then chose \[g=(u_{[-n-r,-n-1]}, u_{[-n,-1]},
u_{[0,p-1]}).\] If on the other hand $r$ is zero and $p$ nonzero we
choose \[g=(u_{[-n-s-1,-n-2]}, u_{[-n-1,-1]}, u_{[0,p-1]}),\] where
$s = |\tau(u_{[-n-1]})| -1$ which is nonzero. If finally $p$ is zero
and $r$ nonzero we choose \[g=(u_{[-n-r,-n-1]}, u_{[-n,0]},
u_{[1,s]}),\] where $s = |\tau(u_{[0]})| -1$ which is nonzero as
well. Note that due to primitivity, we cannot have both $p$ and $r$
zero. The theorem now follows in each case from iterating lemma
\ref{adopt}, making use of the fact that the rl graph is subfixed in
the second case and that the ll graph is subfixed in the third case.
\hspace{\stretch{1}}$\Box$

Let us shortly consider the usefulness of this result: Given a
substitution it is often easy to find some special sequences using
generators, e.g., any two generators with identical right wings but
disagreeing letters in the center complete to left special sequences
modulo orbit equivalence. On the other hand, this result tells us
that under certain circumstances all special sequences can be
obtained in this way. And since the results from the previous
section gives us some measure of control over the generators, we are
now in a better position to face the special sequences of a
substitution. One possible application could be to calculate special
sequences and thereby configuration graphs for arbitrary
substitutions, but this is already done very well in \cite{eilers},
indeed the present section steals heavily from this source. Instead
we shall use our results to produce certain substitutions with
desirable properties such as having a particular configuration
graph; this is the object of the next section.

\section{The Zorro Algorithm}
\subsection{Miscellaneous Tools}
This subsection contains miscellaneous minor results that are needed
in the proof the Zorro Algorithm. While the results are (probably)
true, they may appear unmotivated and rather out of context. Do not
worry though, all will be clear in due time.

\begin{lem} \label{aper}
Let $\tau$ be any primitive substitution. If $\tau$ has either a
left or a right special sequence then it is aperiodic.
\end{lem}
\emph{Proof:} Suppose it has a left special sequence, this provides
us with sequences $x, y \in X_\tau$ with $x_{[-1]} \neq y_{[-1]}$
and $x_{[0,\infty[} = y_{[0,\infty[}$. Assume now that $\tau$ is
periodic, this implies that $x$ and $y$ are each periodic, let $n,
m$ be the lengths of their periods. But then both sequences are
periodic with periods of length $nm$ as well which is an obvious
contradiction. \hspace{\stretch{1}}$\Box$

\begin{prop}
Let $\tau$ be any substitution with subfixed ll and rl graphs. We
have that
\[\mathcal{L}_2(\tau) = \left\{u \in \mathcal{A}_2\mid \exists a \in
\mathcal{A} : u \dashv \tau(a) \right\} \cup \left\{u \in
\mathcal{A}_2 \mid \exists a \in \mathcal{A} : u \dashv \tau^2(a)
\right\}.\]
\end{prop}
\emph{Proof:} Any member of the right hand side is a member of the
left hand side by definition. Now let $u \in \mathcal{L}_2(\tau)$,
by definition we have $a \in \mathcal{A}$ and $n \in \mathbb{N}$
with $u \dashv \tau^n(a)$ and we may chose $a$ and $n$ such that $n$
is minimal. Assume for the sake of contradiction that $n \geq 3$.
This implies that there can be no letter $b \dashv \tau^{n-1}(a)$
with $u \dashv \tau(b)$, nor any letter $b \dashv \tau^{n-2}(a)$
with $u \dashv \tau^2(b)$. But this again implies that there exist
$v,w \in \mathcal{A}^+$ with $vw = \tau^{n-2}(a)$ and with $u =
\mathrm{rl}(\tau^2(v))\mathrm{ll}(\tau^2(w))$. But since the ll and
rl graphs are subfixed we have that
\[\mathrm{rl}(\tau^2(v))\mathrm{ll}(\tau^2(w)) =
\mathrm{rl}(\tau(v))\mathrm{ll}(\tau(w)),\] which implies the
contradiction $u \dashv \tau^{n-1}(a)$. \hspace{\stretch{1}}$\Box$

This result can be generalized to word lengths higher than 2. We
are, however, more interested in the following corollary:
\begin{cor} \label{l2}
Let $\tau$ be any substitution with subfixed ll and rl graphs. Let
\[W = \left\{u \in \mathcal{A}_2\mid \exists a \in \mathcal{A} : u
\dashv \tau(a) \right\}.\] We have that \[\mathcal{L}_2(\tau) = W
\cup \left\{\mathrm{rl}(\tau(a))\mathrm{ll}(\tau(b)) \mid ab \in
W\right\}.\]
\end{cor}

\subsection{The Theorem and the Algorithm}
\begin{defn}
A bipartite graph is said to be \emph{undecided} if it has the
following properties:
\begin{itemize}
\item[(i)] There are no lonely vertices, i.e., any vertex has one or more outgoing edges.
\item[(ii)] There are no lonely edges, i.e., for any edge there exists another edge with one or both vertices mutual.
\item[(iii)] There exists a left vertex with two outgoing edges.
\item[(iv)] There exists a right vertex with two outgoing edges.
\end{itemize}
\end{defn}
We say that a primitive aperiodic substitution \emph{realizes} its
configuration graph and in general that a bipartite graph is
\emph{realizable} if there exists a primitive, aperiodic
substitution realizing it. The following theorem is the conclusion
to much of our work:
\begin{thm}
A bipartite graph is realizable if and only if it is undecided.
Indeed, for any bipartite undecided graph the Zorro algorithm
described below will compute a primitive, aperiodic substitution
realizing it.
\end{thm}
\emph{Proof:} Note initially that by the definitions and results of
subsection \ref{specials} it is immediate that any realizable graph
is undecided. To prove the other way round, we shall first state the
Zorro algorithm with a few examples and then afterwards consider
that it actually produces the desired substitutions.

Consider the following three bipartite graphs:
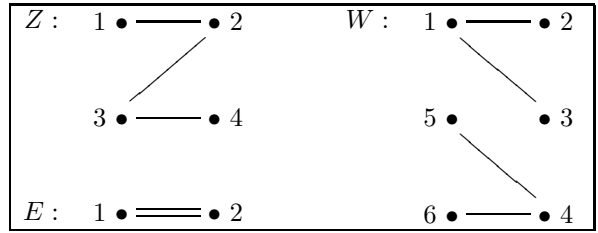
\[\begin{array}{|rlrl|} \hline

Z:&\xymatrix{1\ \bullet \ar@{-}[r] & \bullet\ 2 \\
3\ \bullet \ar@{-}[ur] \ar@{-}[r] & \bullet\ 4} \hspace{6ex}

 &  W: & \xymatrix{1\ \bullet \ar@{-}[r] \ar@{-}[dr] &
\bullet\ 2 \\ 5\ \bullet \ar@{-}[dr] & \bullet\ 3 \\6\ \bullet
\ar@{-}[r] & \bullet\  4} \vspace{-3ex}
\\  E:
&\xymatrix{1\ \bullet \ar@{-}@<+0.3ex>[r] \ar@{-}@<-0.3ex>[r] &
\bullet\ 2} \\ \hline
\end{array}\]
Now let $G$ be any bipartite undecided graph. It follows from parts
(iii) and (iv) of the definition that $G$ must contain one or more
of the above graphs as a subgraph. The algorithm has three cases
corresponding to these three subgraphs, each of these cases proceeds
according to the following common recipe but with slightly differing
ingredients\footnote{Incidentally, the algorithm is named after the
particular shape of the $Z$ graph, this was the first case solved.}:
\begin{enumerate}
\item The first part simply states an \emph{initial substitution}
that realizes the given subgraph. The alphabet has one letter
corresponding to each vertex in the subgraph but also contains
additional letters that do not correspond to vertices. The following
three steps will gradually extend the initial substitution such that
the final result realizes $G$.
\item Remaining vertices are added now: For each vertex in $G$ not
in the subgraph, we add a new letter to our alphabet. The value of
our substitution at these new letters are assigned according to
\emph{left} and \emph{right patterns} for left respectively right
vertices. To be precise, the value of a new letter corresponding to
a left vertex is obtained by postfixing the word produced by the
left pattern with the new letter itself, right letters are treated
symmetrically.
\item Then the first edges: For each pair of vertices
that are presently unconnected but are connected in $G$ we add the
first (possibly only) edge by inserting the two letter word
consisting of the two letters corresponding to the left respectively
right vertex at the \emph{insertion point} specified as part of the
initial substitution.
\item And finally the remaining edges: For any two vertices that are
already connected but lack the number of edges present in $G$, we
add a new letter to our alphabet for each missing edge. The value of
the substitution at such a new letter is obtained by taking first
the value of the substitution at the letter corresponding to the
left vertex minus the rightmost letter, then adding the new letter
and finally the value of the substitution at the letter
corresponding to the right vertex minus the leftmost letter. All new
letters produced in this step are finally added directly as one
letter words at the insertion point.
\end{enumerate}

As a start, let us specify the initial substitution with insertion
point and left and right patterns in the case of the subgraph $Z$,
which is the easiest case:
\[\begin{array}{|lr|} \hline
\begin{array}{rcl}
1 & \mapsto & 2245\boldsymbol{1} \\
2 & \mapsto & \boldsymbol{2}45133 \\
3 & \mapsto & 222451\boldsymbol{3} \\
4 & \mapsto & \boldsymbol{4}51333 \\
5 & \mapsto & 22224\boldsymbol{5}\mid13333 \\
\end{array} &
\begin{array}{l} \mathrm{left\ pattern}:
\underbrace{22\cdots2}_{5,6,7,\ldots}45\\
\mathrm{right\ pattern}:51\underbrace{33\cdots3}_{5,6,7,\ldots}
\end{array} \\  \hline
\end{array}\]
A few words on the notation: The insertion point is specified by a
vertical line, in this case in the middle of the value of 5. As a
theoretical convenience we have highlighted letters in values
letters that are identical to the source letter, this is of no
importance when applying the algorithm. The patterns produce words
of increasing length, i.e., the first word produced by the left
pattern in this case is 2222245, the next 22222245 and so on.
Finally note that the letters 1 though 4 corresponds to the vertices
of $Z$ whereas the letter 5 does not correspond to any vertex.

An example is due, indeed we should very much like to realize the
following graph:
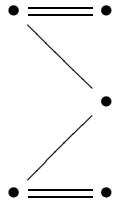
\[\xymatrix{\bullet \ar@{-}@<+0.3ex>[r] \ar@{-}@<-0.3ex>[r]
\ar@{-}[dr] & \bullet \\ & \bullet\\ \bullet \ar@{-}@<+0.3ex>[r]
\ar@{-}@<-0.3ex>[r] \ar@{-}[ur] & \bullet} \] Luckily, it is
undecided. Initially we need to identify which of the three graphs
that are contained in this graph. As it happens, both the $Z$ and
$E$ are subgraphs. For didactic reasons we chose to carry on with
$Z$, but choosing $E$ would have produced a realizing substitution
as well. But then we have an initial substitution and step one of
the algorithm is complete and leaves us with the following
substitution and its configuration graph:
\[
\begin{array}{rcl}
1 & \mapsto & 2245\boldsymbol{1} \\
2 & \mapsto & \boldsymbol{2}45133 \\
3 & \mapsto & 222451\boldsymbol{3} \\
4 & \mapsto & \boldsymbol{4}51333 \\
5 & \mapsto & 22224\boldsymbol{5}\mid13333 \\ \vspace{-14ex}
\end{array}\hspace{6ex}
\xymatrix{1\ \bullet \ar@{-}[dr] & \\ & \bullet\ 2\\
3\ \bullet \ar@{-}[r] \ar@{-}[ur] & \bullet\ 4 }
\]
Notice a two things here: The highlighted symbols and the insertion
point in the substitution are of course not a part of the
substitution but rather theoretically convenient layout, just as the
letters labeling the vertices. Also notice that $Z$ does not occur
as subgraph of our graph in an unambiguous way, indeed we could have
chosen to let the vertices of $Z$ coincide with all vertices except
the lower right instead. This, like the choice between $Z$ and $E$
at step 1, does not matter, all choices will produce realizing, if
not necessarily identical, substitutions. As for step two, we need
to introduce one more vertex, this is done by adding the letter 6 to
our alphabet and assigning it the value 65133333 in accordance with
the right pattern since it is a right vertex. We now have the
following substitution and corresponding graph as conclusion to step
2:
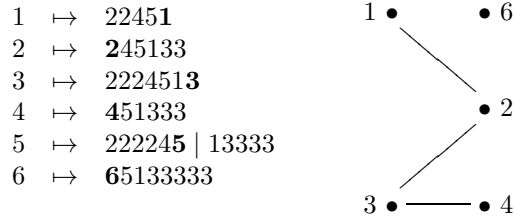
\[
\begin{array}{rcl}
1 & \mapsto & 2245\boldsymbol{1} \\
2 & \mapsto & \boldsymbol{2}45133 \\
3 & \mapsto & 222451\boldsymbol{3} \\
4 & \mapsto & \boldsymbol{4}51333 \\
5 & \mapsto & 22224\boldsymbol{5}\mid13333 \\
6 & \mapsto & \boldsymbol{6}5133333 \\ \vspace{-18ex}
\end{array}\hspace{6ex}
\xymatrix{1\ \bullet \ar@{-}[dr] & \bullet\ 6 \\ & \bullet\ 2\\
3\ \bullet \ar@{-}[r] \ar@{-}[ur] & \bullet\ 4 }
\]
Notice about this step that while the substitution above corresponds
to the graph in algorithmic terms it does not realize it. This is a
slight inconvenience that applies to step two only, essentially it
is caused by adding lonely vertices to the original graph and
thereby wrecking havoc upon its undecidability. As for step three,
we need to add just one edge between the vertices 1 and 6. This is
easily done by adding the two letter word 16 at the insertion point:
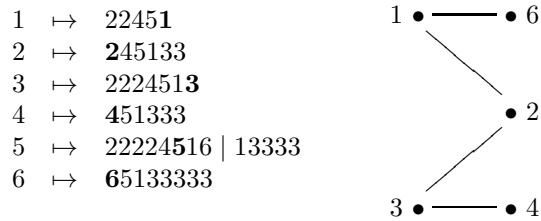
\[
\begin{array}{rcl}
1 & \mapsto & 2245\boldsymbol{1} \\
2 & \mapsto & \boldsymbol{2}45133 \\
3 & \mapsto & 222451\boldsymbol{3} \\
4 & \mapsto & \boldsymbol{4}51333 \\
5 & \mapsto & 22224\boldsymbol{5}16\mid13333 \\
6 & \mapsto & \boldsymbol{6}5133333 \\ \vspace{-18ex}
\end{array}\hspace{6ex}
\xymatrix{1\ \bullet \ar@{-}[r] \ar@{-}[dr] & \bullet\ 6 \\ & \bullet\ 2\\
3\ \bullet \ar@{-}[r] \ar@{-}[ur] & \bullet\ 4 }
\]

Finally, we need to add two more edges between already connected
vertices: One more between vertices 1 and 6 and the final between
the vertices 3 and 4. The first is added by introducing the new
letter 7 and assigning it the value 2245 followed by 7 itself
followed by 5133333, i.e., the unlikely long value of 224575133333.
Similarly the final edge is added by introducing the letter 8 and
assigning it the value 222451851333. Both these two new letters are
added at the insertion point and the fourth and final step of the
algorithm is complete:
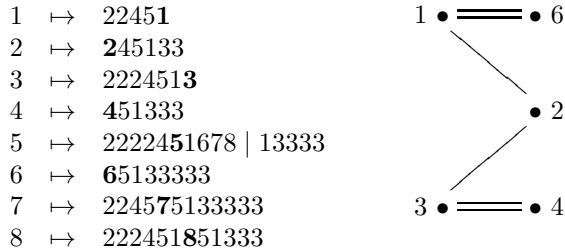
\[
\begin{array}{rcl}
1 & \mapsto & 2245\boldsymbol{1} \\
2 & \mapsto & \boldsymbol{2}45133 \\
3 & \mapsto & 222451\boldsymbol{3} \\
4 & \mapsto & \boldsymbol{4}51333 \\
5 & \mapsto & 22224\boldsymbol{5}1678\mid13333 \\
6 & \mapsto & \boldsymbol{6}5133333 \\
7 & \mapsto & 2245\boldsymbol{7}5133333 \\
8 & \mapsto & 222451\boldsymbol{8}51333 \\ \vspace{-23ex}
\end{array}\hspace{6ex}
\xymatrix{1\ \bullet \ar@{-}@<+0.3ex>[r] \ar@{-}@<-0.3ex>[r] \ar@{-}[dr] & \bullet\ 6 \\ & \bullet\ 2\\
3\ \bullet \ar@{-}@<+0.3ex>[r] \ar@{-}@<-0.3ex>[r] \ar@{-}[ur] &
\bullet\ 4 }
\] \vspace{1ex}

The example concluded, let us now state the initial substitution
etc. for the remaining two cases. First the case of the subgraph
$W$:
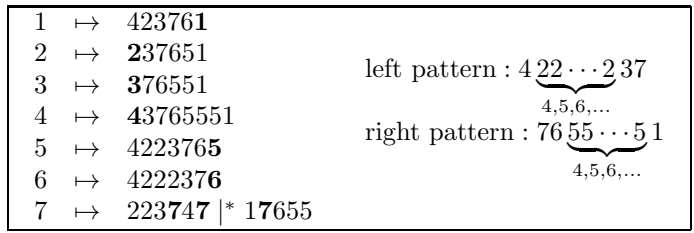
\[\begin{array}{|lr|} \hline
\begin{array}{rcl}
1 & \mapsto & 42376\boldsymbol{1} \\
2 & \mapsto & \boldsymbol{2}37651 \\
3 & \mapsto & \boldsymbol{3}76551 \\
4 & \mapsto & \boldsymbol{4}3765551 \\
5 & \mapsto & 422376\boldsymbol{5} \\
6 & \mapsto & 422237\boldsymbol{6} \\
7 & \mapsto & 223\boldsymbol{7}4\boldsymbol{7}\mid^*1\boldsymbol{7}655 \\
\end{array} &
\begin{array}{l} \mathrm{left\ pattern}:
4\underbrace{22\cdots2}_{4,5,6,\ldots}37\\
\mathrm{right\ pattern}:76\underbrace{55\cdots5}_{4,5,6,\ldots}1
\end{array} \\  \hline
\end{array}\]
As hinted by the star next to the insertion point, there is one
peculiarity to this case as compared to the two others: All words
inserted at the insertion point, whether at step three or four in
the algorithm, need to be followed by the letter 7, e.g., if the
algorithm tells us to insert the words 53, 8 and 9 at the insertion
point, then we need to insert
$53\boldsymbol{7}8\boldsymbol{7}9\boldsymbol{7}$ and not just their
concatenation 5389 as we would in the other two cases. This is
caused, in a sense, by the graph $W$ being disconnected, the symbol
7 works as bridge between the parts. The final case of the subgraph
$E$ completes the definition of the algorithm:
\[\begin{array}{|lr|} \hline
\begin{array}{rcl}
1 & \mapsto & 253425\boldsymbol{1} \\
2 & \mapsto & \boldsymbol{2}513451 \\
3 & \mapsto & 25\boldsymbol{3}425\boldsymbol{3} \\
4 & \mapsto & \boldsymbol{4}513\boldsymbol{4}51 \\
5 & \mapsto & 2\boldsymbol{5}1134 \mid 3422\boldsymbol{5}1 \\
\end{array} &
\begin{array}{l} \mathrm{left\ pattern}:
25\underbrace{11\cdots1}_{3,4,5,\ldots}34\\
\mathrm{right\ pattern}:34\underbrace{22\cdots2}_{3,4,5,\ldots}51
\end{array} \\  \hline
\end{array}\]

As conclusion to our description of the algorithm we provide two
more examples, one for each of the graphs $W$ and $E$. We shall not
go into the same level of detail as before, but rather just present
the desired graphs and then state the results of running the
algorithm. The two undecided graphs we would like to realize are:
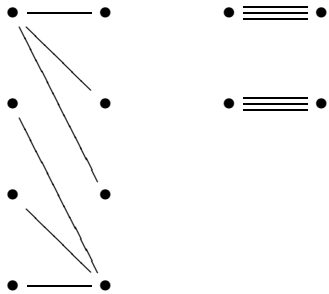
\[\begin{array}{lr}
\xymatrix{\bullet \ar@{-}[r], \ar@{-}[dr], \ar@{-}[ddr] & \bullet \\
          \bullet \ar@{-}[ddr] & \bullet \\
          \bullet \ar@{-}[dr] & \bullet \\
          \bullet \ar@{-}[r] & \bullet} &
\hspace{6ex}\xymatrix{\bullet \ar@{-}@<-0.5ex>[r] \ar@{-}[r]
\ar@{-}@<+0.5ex>[r] & \bullet \\ \bullet \ar@{-}@<-0.5ex>[r]
\ar@{-}[r] \ar@{-}@<+0.5ex>[r] & \bullet}
\end{array}\]

The first contains the graph $W$ and running the algorithm gives us
the following result:
\[
\begin{array}{rcl}
1 & \mapsto & 42376\boldsymbol{1} \\
2 & \mapsto & \boldsymbol{2}37651 \\
3 & \mapsto & \boldsymbol{3}76551 \\
4 & \mapsto & \boldsymbol{4}3765551 \\
5 & \mapsto & 422376\boldsymbol{5} \\
6 & \mapsto & 422237\boldsymbol{6} \\
7 & \mapsto & 223\boldsymbol{7}4\boldsymbol{7}18\boldsymbol{7}94\boldsymbol{7}\mid^*1\boldsymbol{7}655 \\
8 & \mapsto & \boldsymbol{8}7655551 \\
9 & \mapsto & 4222237\boldsymbol{9} \\ \vspace{-26ex}
\end{array}\hspace{6ex}
\xymatrix{1\ \bullet \ar@{-}[r], \ar@{-}[dr], \ar@{-}[ddr] & \bullet\ 2 \\
          9\ \bullet \ar@{-}[ddr] & \bullet\ 3 \\
          5\ \bullet \ar@{-}[dr] & \bullet\ 8 \\
          6\ \bullet \ar@{-}[r] & \bullet\ 4}
\]
Notice here how the words 18 and 94 are followed by the letter 7 as
specified above. The final example gives the following
result:\newpage
\[
\begin{array}{rcl}
1 & \mapsto & 253425\boldsymbol{1} \\
2 & \mapsto & \boldsymbol{2}513451 \\
3 & \mapsto & 25\boldsymbol{3}425\boldsymbol{3} \\
4 & \mapsto & \boldsymbol{4}513\boldsymbol{4}51 \\
5 & \mapsto & 2\boldsymbol{5}1134 67890\mid 3422\boldsymbol{5}1 \\
6 & \mapsto & 2511134\boldsymbol{6} \\
7 & \mapsto & \boldsymbol{7}3422251 \\
8 & \mapsto & 253425\boldsymbol{8}513451 \\
9 & \mapsto & 2511134\boldsymbol{9}3422251 \\
0 & \mapsto & 2511134\boldsymbol{0}3422251 \\ \vspace{-28ex}
\end{array}\hspace{6ex}
\xymatrix{1\ \bullet \ar@{-}@<-0.5ex>[r] \ar@{-}[r]
\ar@{-}@<+0.5ex>[r] & \bullet\ 2 \\ 6\ \bullet \ar@{-}@<-0.5ex>[r]
\ar@{-}[r] \ar@{-}@<+0.5ex>[r] & \bullet\ 7}
\] \vspace{15ex}

Having thus stated and exemplified the algorithm it is time to prove
that it does indeed produce a primitive, aperiodic substitution
realizing a given graph. We shall not go through all painstaking
details three times. Instead, we list the properties that need to be
verified for all cases and for each of these properties describe the
general strategy used to verify it. And we shall, of course, verify
a few of these properties in full detail for some of the cases.

The first issue to consider is that of primitivity. This is
fundamental to all our workings and luckily it holds easily for all
substitutions since they all contain a particular letter with the
property that its value contains the entire alphabet and it is
itself contained in the value of all letters. This is the letter
with the insertion point. Notice that this property also holds after
adding additional letters according to step 2 since these are all
added at the insertion point in part 3 by part (i) of the definition
of an undecided graph. The next basic issue is aperiodicity, but
this is easily handled by lemma \ref{aper} since left or right
special sequences are easily constructed from generators in all the
initial substitutions. As an example, the generators
$(253425,12,513451)$ and $(253425,34,513451)$ from case $E$ provide
us with both left and right special sequences.

Having dealt with the basics, we now check that the produced
substitutions are prefix as well as postfix free, this implies that
they are segregating with least left and right segregating numbers
both 1. With this in mind, we furthermore verify that all the four
graphs ll, rl, ls and rs are subfixed. This is where the weird
patterns used in step 2 are justified since they oversee that these
properties, that hold for the initial substitutions, are maintained
through the steps 2, 3 and 4 of the algorithm. Take as an example
the case $W$: The initial substitution is easily prefix and postfix
free and some checking shows that the four graphs are all subfixed.
Now let us add a left vertex as an example of the effects of step 2,
we get $8 \mapsto 42222378$. On the right hand side the new unique
letter 8 protects from trouble. And the left pattern ensures not
only that the substitution remains prefix free but also that the ll
and in particular the ls graph remain subfixed. Step 3 changes
nothing and the letters introduced in step 4 also ends up being
compatible with the state of affairs. Taking some time to verify
these things also gives some idea of why the produced substitutions
tend to be lengthy.

We now have primitive, aperiodic, segregating substitutions with the
four graphs ll, rl, ls and rs subfixed. And indeed, we are going
strong, these are exactly the prerequisites of theorem \ref{exists}.
The next consideration is to identify the set of basic generators
for each substitution and from these verify that the desired graph
is actually realized. Let us consider an example to simplify things:
Letting $\tau$ be the initial substitution in the case $Z$ we easily
get by corollary \ref{l2} that $12, 32, 34 \in \mathcal{L}_2(\tau)$
but $14 \notin \mathcal{L}_2(\tau)$, which again easily gives us the
following basic generators:
\[ (22224,5,13333),\ (2245,12,45133),\ (222451,32,45133),\
(222451,34,51333). \] This immediately implies that the orbit
classes containing the completions of the three last generators are
special and by proposition \ref{simsim} different. Furthermore, by
theorem \ref{exists} and corollary \ref{simrid} these are the only
special orbit classes. And by corollary \ref{simrid} the completions
of the second and third are right tail equivalent whereas the
completion of the fourth isn't right tail equivalent with any of the
others; similarly the completions of the third and fourth are left
tail equivalent but the second is excluded. Summing up, we have
proved that the initial substitution actually does realize the $Z$
graph, and in general that, because of our careful preparations
above, the configuration graph is easily read off from the set of
basic generators.

The general idea is now that any left vertex corresponds to a letter
with a value consisting of a unique left part not containing the
letter itself followed by the letter. This correspondence is set up
in the initial substitution and is maintained through step 2 by the
left pattern. The situation is symmetrical for the right vertices.
Step 2 does thus not in itself produce any new generators, since the
centers of the potential generators are not in the language yet.
This setup makes the adding of vertices at step 3 very easy though,
just extend the language by adding words at the insertion point,
only we have to take some care in the case of case $W$ not to
introduce unwanted generators. Note that by part (ii) of the
definition of an undecided graph we are ensured that all edges share
a vertex with some other edge, this ensures that the generators we
add in this step become special and thus actually figure in the
graph. At step 4 we want to add an additional edge between already
connected vertices, this is easily done by introducing a new
generator with left and right wings corresponding to the vertices
but with a new center, and remembering to add it to the language. To
satisfactorily verify the algorithm one of course needs to check
very carefully that no unwanted two letter words enter the language
during the steps 2 through 4, since this would give an undesired
edge, we shall refrain from doing this in writing.
 \hspace{\stretch{1}}$\Box$


\begin{thebibliography}{9}
\bibitem{eilers} T. M. Carlsen and S. Eilers: \emph{A Graph Approach
to Computing Nondeterminacy in Substitutional Dynamical Systems},
submitted, \texttt{www.math.ku.dk/\~{}eilers/papers/cei.html}, 2002.
\bibitem{applet} T. M. Carlsen and S. Eilers: \emph{Java applet},
\texttt{www.math.ku.dk/\~{}eilers/papers/cei.html}, 2002.
\bibitem{salomaa} G. Rozenberg and A. Salomaa: \emph{Mathematical
Theory of L systems}, Academic Press Inc., 1980.
\bibitem{durand} F. Durand, B. Host and C. Skau:
\emph{Substitutional Dynamical Systems, Bratelli Diagrams and
Dimension Groups}, Ergodic Theory Dynam. Systems $\boldsymbol{19}$
(1999), no. 4, 953-993.
\end{thebibliography}
\end{document}